\documentclass[subscriptcorrection,upint,varvw,barcolor=Goldenrod3,balance,hyphenate,pdf-a,nolists]{asmejour}

\usepackage{cleveref}
\usepackage{booktabs,xcolor,siunitx}
\usepackage{mathtools,algorithm}
\usepackage[noend]{algpseudocode}
\usepackage{algorithm}

\newcommand{\overbar}[1]{\mkern 1.5mu\overline{\mkern-1.5mu#1\mkern-1.5mu}\mkern 1.5mu}
\newcommand{\f}[2]{\frac{#1}{#2}}
\newcommand{\mb}[1]{\mathbf{#1}}

\DeclareMathAlphabet\mathbfcal{OMS}{cmsy}{b}{n}
\renewcommand{\d}{\mathop{}\!\mathrm{d}} 
\newcommand{\p}{\partial}

\hypersetup{%
pdfauthor={Sicheng He, Shugo Kaneko, Max Howell, Nan Li, Joaquim R. R. A. Martins},                       		 
pdftitle={MDO Control}, 
pdfkeywords={MDO multidisciplinary design optimization LQR cartpole control codesign CCD},
pdfsubject = {Control Co-Design with MDO},
}

\JourName{Mechanical Design}

\begin{document}

\SetAuthorBlock{Sicheng He}{%
Assistant Professor\\
Dept. of Mech., Aero., Biomed. Eng.\\
Univ. of Tennessee, Knoxville, TN 37996 \\
email: sicheng@utk.edu
}
\SetAuthorBlock{Shugo Kaneko}{%
Postdoctoral Fellow\\
Dept. of Aerospace Eng. \\
Univ. of Michigan, Ann Arbor, MI 48109 \\
email: shugok@umich.edu
}
\SetAuthorBlock{Max Howell}{%
Graduate Student\\
Dept. of Mech., Aero., Biomed. Eng.\\
Univ. of Tennessee, Knoxville, TN 37996 \\
email: mhowel30@vols.utk.edu
}
\SetAuthorBlock{Nan Li}{%
Assistant Professor\\
Dept. of Aero.Eng. \\
Auburn University, Auburn, AL 36849 \\
email: nanli@au.edu
}
\SetAuthorBlock{Joaquim R. R. A. Martins}{%
Pauline M. Sherman Professor\\
Dept. of Aerospace Eng.\\
Univ. of Michigan, Ann Arbor, MI 48109 \\
email: jrram@umich.edu
}

\title{Efficient Adjoint-based Design Optimization with Optimal Control}


\begin{abstract}
Multidisciplinary engineering system design typically employs a sequential process, progressing from system dynamics to design variables and control. 
However, this process is inefficient and may lead to a suboptimal design. 
We propose formulating the optimal control and multidisciplinary design optimization (MDO) problems as a single problem with linear quadratic regulator (LQR) control.
We use the coupled adjoint method to compute the design variable derivatives, which are critical for gradient-based design optimization.
The computational cost of the derivative computation using the adjoint method is independent of the number of design variables, making it suitable for large-scale problems. 
We show that the coupled adjoint can be solved indirectly and more efficiently by solving three smaller adjoint equations that leverage the feedforward structure of the problem.
We demonstrate this new approach on two test problems: design optimization of a classic cart-pole problem and the aerodynamic shape of a quadrotor blade.
For the quadrotor blade design problem, we reduce the control cost by 10\% by optimizing the blade for a specific control task with a slight penalty in steady hovering power consumption.
\end{abstract}

\date{Version \versionno, \today}
\maketitle

\section{Introduction}

Multidisciplinary engineering system design typically employs a sequential process, where engineers first optimize the physical system design variables (e.g., structural dimensions, material properties) and then design a controller for the fixed plant. 
Engineers iterate this process several times to converge on an optimized design.
However, sequential optimization typically yields suboptimal designs~\cite[Sec.~13.1]{Martins2022}. 
Control co-design (CCD) addresses this issue.
CCD is a type of multidisciplinary design optimization (MDO) problem that optimizes the system and control design simultaneously.
CCD has been applied to structural design~\cite{Eastep1987a,RAO1988a,HAFTKA1990a,Belvin1990a,Sunar1993a}, automobiles~\cite{Fathy2003a,Allison2014a}, aircraft design~\cite{Gupta2020a,Fasse2022a,Kaneko2024b}, wind turbines~\cite{Deshmukh2015a,Pao2021a,LpezMuro2022a}, chaotic dynamical systems~\cite{Cagliari2021a}, and robotics~\cite{Asada1991a,Fadini2021a}.
Allison and Herber~\cite{Allison2014b} and Garcia-Sanz~\cite{GarciaSanz2019a} present broader reviews of CCD.

Open and closed-loop control strategies have both been used in complex CCD problems. 
Open-loop co-design explicitly optimizes the control trajectory, while closed-loop control finds an optimal control law that is applied to the system using state feedback. 
Open-loop CCD has been applied to spacecraft~\cite{Hale1985}, automobiles~\cite{Allison2014a}, and many other systems. 
However, systems that employ open-loop control often lack robustness and may fail due to external disturbances. 
Closed-loop control strategies, such as the linear quadratic regulator (LQR), are advantageous because they can self-correct using state feedback.
Closed-loop CCD has been used to optimize controlled space structures with LQR~\cite{Belvin1990a}. 
Aksland et al.~\cite{Aksland2023} outline an approach to robust closed-loop CCD using both stochastic and deterministic optimization.

CCD optimization problems can be solved using multiple approaches, depending on how we partition and couple the design and control problems~\cite{Reyer2001a}.
Two of the most common approaches are the \emph{simultaneous} formulation and the \emph{nested} formulation~\cite{Herber2018a, Kaneko2025a}.
The simultaneous formulation solves a monolithic optimization problem where one optimizer drives both design and control variables simultaneously.
The nested formulation divides the co-design problem into two levels: the inner loop computes the optimal control for a given plant design, while the outer loop optimizes the plant variables~\cite{Fathy2003a}. 
Herber and Allison~\cite{Herber2018a} and Sundarrajan and Herber~\cite{Sundarrajan2021a} review and compare the simultaneous and nested approaches.

The simultaneous formulation is often more efficient for open-loop CCD problems~\cite{Kaneko2025a}, although it is also applicable to closed-loop problems.
Allison et al.~\cite{Allison2014a} developed the direct transcript (DT) method, which imposes the dynamics (equations of motion) as a series of equality constraints.
Hence, the DT method allows the optimizer to find a physical solution of the dynamical system while simultaneously optimizing the design.
DT corresponds to the \emph{all-at-once} (AAO) or \emph{simultaneous analysis and design} (SAND) architecture in the MDO architecture terminology~\cite{Martins2013}. 

One disadvantage of DT-based formulation is that the optimization problem size can become prohibitively large and expensive to solve when there are many state and control variables, or when considering a long time horizon with a small time step.
This is because DT treats both the state and control variables at each time step as optimization variables. 
The simultaneous formulation can also use the multidisciplinary feasible (MDF) architecture~\cite{Martins2013} by using a time marching solver or a nonlinear solver to simulate dynamics.
However, the DT-based approach tends to be more robust and computationally efficient than the MDF-based approach in simultaneous formulation~\cite{Kaneko2025a}.

In the nested formulation, the problem size of outer-loop design optimization is significantly smaller than that of the simultaneous formulation because it only optimizes the plant design variables.
However, it now has to solve inner-loop optimal control at every outer-loop optimization iteration.
Therefore, the inner-loop optimal control problem must be solved computationally efficiently for the nested formulation to outperform the simultaneous formulation~\cite{Herber2018a}.
This is possible when the inner-loop problem has certain characteristics.
For example, LQR control provides an analytical optimum for the control problem, in which case an iterative optimizer is not needed for the inner loop~\cite{Fathy2001a}. 
One of the main challenges in the nested formulation is the computation of derivatives for the optimal control, which is required by the outer-loop optimization when using a gradient-based method.
Without the derivatives of the optimal control, we must use gradient-free optimization for the outer loop, which cannot solve large-scale optimization problems.
Eastep et al.~\cite{Eastep1987a} solved a structural optimization problem using the nested formulation, where the inner level is the optimal control problem solved with LQR and the outer level is an eigenvalue problem.
The eigenvalue derivative with respect to the design variable is computed using a forward method, resulting in a computational cost that is proportional to the number of design variables.
For more advanced cases, it may be necessary to solve equations to locate the equilibrium point.
More recent studies used finite difference approximations to compute the derivatives of optimal control~\cite{Herber2018a, Sundarrajan2021a}.
However, finite differences are not efficient when the number of plant design variables is large~\cite[Sec.~6.4]{Martins2022}.
Cunis et al.~\cite{Cunis2023} discussed the analytical derivative computations of open-loop optimal control.

In this paper, we develop an efficient derivative computation method using the adjoint method for the nested closed-loop CCD formulation to enable gradient-based optimization with a large number of design variables. 
Three sets of coupled equations are considered: equilibrium point equations, optimal condition equations and closed-loop unsteady equations.
The optimal condition equations are modeled as an algebraic Riccati equation (ARE), which comes from the LQR used in this study. 
Solving the coupled adjoint equation for this problem is challenging due to its large size, and a specialized solver needs to be developed (e.g., coupled Newton--Krylov with specially designed preconditioners).

We exploit the feed-forward structure of the coupled adjoint equations by solving three individual adjoint equations. 
This reduces the computational cost and removes the need for a new coupled adjoint solver, reducing implementation effort. 
The proposed method can be used to solve CCD problems with a large number of design variables. 
However, in this paper, we did not address the possible high-dimensionality of the state variable. 
Although it can be handled efficiently for the equilibrium point equations and the closed-loop unsteady equations, solving large-scale ARE and its corresponding adjoint equation, which is often posed as the Lyapunov equation, remains an open research problem.

There are several relevant works on solving large-scale AREs and Lyapunov equations, including the matrix equation sparse solver (M.E.S.S.)~\cite{mess2025}. 
Benner et al.~\cite{benner2019efficient,Benner2020} address the computational challenges of solving large-scale AREs arising in LQR problems for systems modeled by index-2 differential-algebraic equations (DAEs), such as those encountered in flow control and constrained fluid-structure interactions. 
The central contribution is an inexact low-rank Newton-ADI method, which combines Newton iterations with the alternating direction implicit (ADI) method while maintaining low-rank representations throughout the process.
A key strength of the method is its ability to handle extremely large problem sizes---with state-space dimensions exceeding $10^5$---without forming or storing full system matrices. 
This is made possible by exploiting the sparsity and low-rank structure of the linearized operators. 

The paper is organized as follows.
In \cref{sec:control}, we present the CCD problem formulation with all its components.
Then, we develop the efficient adjoint formulas for the derivative computation in \cref{sec:oc-derivatives}.
We showcase the performance of the proposed algorithm in \cref{sec:numerical_example}.
Finally, we conclude the paper in \cref{sec:conclusion}.

\section{Control co-design (CCD) problem formulation}
\label{sec:control}

We consider a parametrized dynamical system, $\mb{r}_{\mathrm{nl}}(\mb{x}, \mb{d})$, with full-state feedback gain $\mb{G}(\mb{x}_{\mathrm{tgt}})$ acting on control input $\mb{u}$.
Given mechanical design variables $\mb{d}$ and initial state $\mb{x}_{\mathrm{init}}$, the classic optimal control problem minimizes control cost to reach a prescribed target $\mb{x}_{\mathrm{tgt}}$ by varying the full-state feedback matrix $\mb{W}$.
We apply CCD to jointly optimize the mechanical design $\mb{d}$ and control input, adding a second optimization level above the classical problem.
This extended design space improves final design performance.

\subsection{Equilibrium point equation solution}

The target point is an equilibrium point that may or may not require control input to balance external loads.
The equilibrium point is obvious in some cases; however, for more general cases, the detailed states at the equilibrium point are unknown and must be computed.
In this section, we develop equations to compute the state variables and control variables at the equilibrium point.

Suppose the underlying nonlinear dynamic equation with external control input is written as
\begin{equation}
\label{eq:dyn_open}
\dot{\mb{x}} = \mb{r}_{\mathrm{nl}}(\mb{x}, \mb{u}, \mb{d}),
\end{equation}
where $\mb{r}_{\mathrm{nl}}:\mathbb{R}^{n_{\mb{x}}}\times \mathbb{R}^{n_{\mb{u}}} \times \mathbb{R}^{n_{\mb{d}}} \rightarrow \mathbb{R}^{n_{\mb{x}}}$ is the residual form for the nonlinear dynamical system, $\mb{x}\in\mathbb{R}^{n_{\mb{x}}}$ is the state variable with $n_{\mb{x}}$ denotes its dimension, $\mb{u}\in \mathbb{R}^{n_{\mb{u}}}$ is the control variable with $n_{\mb{u}}$ denotes its dimension, and $\mb{d}\in\mathbb{R}^{n_{\mb{d}}}$ is the design variable with $n_{\mb{d}}$ denotes its dimension.
For the equilibrium point, $\dot{\mb{x}} = \mb{0}$.
Thus, from \cref{eq:dyn_open}, the equilibrium point solution satisfies
\begin{equation}
\label{eq:dyn_open_eql}
\mb{r}_{\mathrm{nl}}(\mb{x}_{\mathrm{tgt}}, \mb{u}_{\mathrm{tgt}}, \mb{d}) = \mb{0}, 
\end{equation}
where $\mb{x}_{\mathrm{tgt}}, \mb{u}_{\mathrm{tgt}}$ denote the state variable and control variable at the equilibrium point, respectively.

We may have partial knowledge of state variables, $\mb{x}_{\mathrm{tgt}}$, and the control variable, $\mb{u}_{\mathrm{tgt}}$.
We can decompose the state and the control variables based on whether they are known \emph{a priori}
\begin{equation}
\mb{x}_{\mathrm{tgt}} = 
\begin{bmatrix}
\hat{\mb{x}}_{\mathrm{tgt}} \\
\check{\mb{x}}_{\mathrm{tgt}}
\end{bmatrix}, \quad
\mb{u}_{\mathrm{tgt}} = 
\begin{bmatrix}
\hat{\mb{u}}_{\mathrm{tgt}} \\
\check{\mb{u}}_{\mathrm{tgt}}
\end{bmatrix},
\end{equation}
where $\hat{\mb{x}}_{\mathrm{tgt}}\in \mathbb{R}^{\hat{n}_{\mb{x}}}$ and $\check{\mb{x}}_{\mathrm{tgt}}\in\mathbb{R}^{\check{n}_{\mb{x}}}$ are the known and unknown state variables, respectively ($\check{n}_{\mb{x}} + \hat{n}_{\mb{x}} = n_{\mb{x}}$); 
$\hat{\mb{u}}_{\mathrm{tgt}}\in\mathbb{R}^{\hat{n}_{\mb{u}}}$ and $\check{\mb{u}}_{\mathrm{tgt}}\in\mathbb{R}^{\check{n}_{\mb{u}}}$ are known and unknown control variables, respectively ($\check{n}_{\mb{u}} + \hat{n}_{\mb{u}} = n_{\mb{u}}$).
We assume that the known states appear before the unknown states because we can always reorder the equations to make sure this is the case.

We assume that we can extract a subset of independent equations from \cref{eq:dyn_open_eql} which has the dimension of $\check{n}_{\mb{x}} + \check{n}_{\mb{u}}$.
This set of equations can be written as
\begin{equation}
\label{eq:steady_state_full}
\check{\mb{r}}_{\mathrm{nl}}(\check{\mb{x}}_{\mathrm{tgt}}, \check{\mb{u}}_{\mathrm{tgt}},  \hat{\mb{x}}_{\mathrm{tgt}}, \hat{\mb{u}}_{\mathrm{tgt}}, \mb{d})=\mb{0},
\end{equation}
where the dimension of the residual, $\check{\mb{r}}_{\mathrm{nl}}$, is $\check{n}_{{\mathrm{nl}}} = \check{n}_{\mb{x}} + \check{n}_{\mb{u}}$.
Defining the vector 
\begin{equation}
\boldsymbol{\theta} = 
\begin{bmatrix}
\check{\mb{x}}_{\mathrm{tgt}} \\ 
\check{\mb{u}}_{\mathrm{tgt}},
\end{bmatrix},
\end{equation}
we can write \cref{eq:steady_state_full} compactly as
\begin{equation}
\label{eq:steady_state_full_compact}
\check{\mb{r}}_{\mathrm{nl}}(\boldsymbol{\theta}; \hat{\mb{x}}_{\mathrm{tgt}}, \hat{\mb{u}}_{\mathrm{tgt}}, \mb{d}) = \mb{0}.
\end{equation}
By solving this nonlinear equation, we can obtain the unknown target state, $\check{\mb{x}}_{\mathrm{tgt}}$, and the unknown target control, $\check{\mb{u}}_{\mathrm{tgt}}$.

There are several special cases in the partition of the state and control variables.
One special case is when the steady-state variable is completely unknown, and the control signal is known \emph{a priori}.
In this case, we have
\begin{equation}
\boldsymbol{\theta} = \mb{x}_{\mathrm{tgt}}, \quad \mb{x}_{\mathrm{tgt}} = \check{\mb{x}}_{\mathrm{tgt}}, \quad \mb{u}_{\mathrm{tgt}} = \hat{\mb{u}}_{\mathrm{tgt}}.
\end{equation}
Another special case is when the full steady-state is specified.
In this case, we have
\begin{equation}
\boldsymbol{\theta} = \mb{u}_{\mathrm{tgt}}, \quad \mb{x}_{\mathrm{tgt}} = \hat{\mb{x}}_{\mathrm{tgt}}, \quad \mb{u}_{\mathrm{tgt}} = \check{\mb{u}}_{\mathrm{tgt}}.
\end{equation}
In this paper, we derive our result for the most general case from \cref{eq:steady_state_full}. 

\subsection{Algebraic Riccati equation (ARE) solution}

After computing the target point states by solving \cref{eq:steady_state_full_compact}, we can develop the optimal control module by linearizing at the target point.
We use the LQR for our model, and find the feedback control coefficients by solving the ARE~\cite[Ch.~3]{Skogestad2005}. 

The linearized dynamical system in the vicinity of the equilibrium point---computed earlier by solving \cref{eq:steady_state_full_compact}---can be written as
\begin{equation}
\label{eq:rlin}
\mb{r}_{\mathrm{lin}}(\delta\mb{x}, \delta\mb{u}, \mb{J}_{\mathrm{tgt}}, \mb{G}_{\mathrm{tgt}})=\dot{\delta\mb{x}} - \mb{J}_{\mathrm{tgt}} \delta\mb{x} - \mb{G}_{\mathrm{tgt}}\delta\mb{u}=\mb{0},
\end{equation} 
where $\delta\mb{x}\in \mathbb{R}^{n}$ is a perturbed state variable measuring the deviation from the equilibrium point, $\delta\mb{u}\in\mathbb{R}^{n_{\mb{u}}}$ is perturbed control input to be determined, $\mb{J}(\mb{x}_{\mathrm{tgt}}, \mb{u}_{\mathrm{tgt}},\mb{d}): \mathbb{R}^{n_{\mb{x}}}\times \mathbb{R}^{n_{\mb{u}}} \times \mathbb{R}^{n_{\mb{d}}}\rightarrow\mathbb{R}^{n_{\mb{x}}\times n_{\mb{x}}}$ is the Jacobian matrix, and $\mb{G}(\mb{x}_{\mathrm{tgt}}, \mb{u}_{\mathrm{tgt}},\mb{d}): \mathbb{R}^{n_{\mb{x}}}\times \mathbb{R}^{n_{\mb{u}}} \times \mathbb{R}^{n_{\mb{d}}}\rightarrow\mathbb{R}^{n_{\mb{x}}\times n_{\mb{u}}}$ is the coefficient matrix relating the control input and the control load.
More specifically, we have
\begin{equation}
\begin{aligned}
\mb{J}_{\mathrm{tgt}} &= \mb{J}\left({\mb{x}_{\mathrm{tgt}}, \mb{u}_{\mathrm{tgt}}, \mb{d}}\right) = \f{\p\mb{r}_{\mathrm{nl}}}{\p \mb{x}}\left({\mb{x}_{\mathrm{tgt}}, \mb{u}_{\mathrm{tgt}}, \mb{d}}\right), \\
\mb{G}_{\mathrm{tgt}} &= \mb{G}\left({\mb{x}_{\mathrm{tgt}}, \mb{u}_{\mathrm{tgt}}, \mb{d}}\right) = \f{\p\mb{r}_{\mathrm{nl}}}{\p \mb{u}}\left({\mb{x}_{\mathrm{tgt}}, \mb{u}_{\mathrm{tgt}}, \mb{d}}\right). \\
\end{aligned}
\end{equation}
When applying LQR to \cref{eq:rlin}, $\delta\mb{u}$ is set as a linear function of $\delta\mb{x}$. 
Then, we have
\begin{equation}
\delta \mb{u} = \mb{W}\delta\mb{x}
\end{equation} 
where $\mb{W}\in\mathbb{R}^{n_{\mb{u}}\times n}$ is the feedback control matrix.


Then, the optimum LQR parameter, $\mb{W}$, is
\begin{equation}
\label{eq:fsfbW}
\mb{W} = - \mb{S}^{-1} \mb{G}_{\mathrm{tgt}}^\intercal \mb{P},
\end{equation}
where $\mb{P}\in\mathbb{R}^{n\times n}$ is found by solving the ARE, and the matrix $\mb{S}\in\mathbb{R}^{n_{\mb{u}}\times n_{\mb{u}}}$ is a positive semi-definite matrice denoted as the control cost matrix.
The ARE is defined by
\begin{equation}
\label{eq:res_ARE}
\mb{R}_{\mathrm{ARE}}(\mb{P}; \mb{J}_{\mathrm{tgt}}, \mb{G}_{\mathrm{tgt}}, \mb{Q}, \mb{S}) = \mb{J}_{\mathrm{tgt}}^\intercal \mb{P} + \mb{P} \mb{J}_{\mathrm{tgt}} - \mb{P} \mb{G}_{\mathrm{tgt}} \mb{S}^{-1} \mb{G}_{\mathrm{tgt}}^\intercal\mb{P} + \mb{Q} = \mb{O},
\end{equation}
where $\mb{Q}\in\mathbb{R}^{n_{\mb{x}}\times n_{\mb{x}}}$ is a positive semi-definite matrix denoted as the state cost matrix.
For the ARE, the state variable is $\mb{P}$, and $\mb{J}_{\mathrm{tgt}}, \mb{G}_{\mathrm{tgt}}, \mb{Q}, \mb{S}$ are coefficient matrices.

\subsection{Closed-loop system equation}

After computing the LQR feedback matrix, $\mb{W}$, we finally have a closed-loop system.
Applying LQR to the original dynamical system \cref{eq:dyn_open}, we obtain the closed-loop system
\begin{equation}
\label{eq:unsteady_control_1}
\mb{r}_{\mathrm{cl}}(\delta\mb{x}; \mb{x}_{\mathrm{tgt}}, \mb{u}_{\mathrm{tgt}}, \mb{W}, \mb{d}) = \dot{\delta\mb{x}} - \mb{r}_{\mathrm{nl}}(\mb{x}_{\mathrm{tgt}} + \delta\mb{x}, \mb{u}_{\mathrm{tgt}} + \mb{W}\delta\mb{x}, \mb{d}) = \mb{0}. 
\end{equation} 
For simplicity, we apply the first-order forward Euler time integration scheme to simulate the dynamical system.
Assuming a finite time-horizon with $n_t$ steps, we have 
\begin{equation}
\label{eq:unsteady_control_2}
\mb{r}_{\mathrm{cl}}^{(i)} = \f{{\delta\mb{x}}^{(i)}- {\delta\mb{x}}^{(i-1)}}{\Delta t} - \mb{r}_{\mathrm{nl}}(\mb{x}_{\mathrm{tgt}} + \delta\mb{x}^{(i)}, \mb{u}_{\mathrm{tgt}} + \mb{W}\delta\mb{x}^{(i)}, \mb{d}) = \mb{0},
\end{equation}
where $i=1, \ldots , n_t$, $\Delta t$ is the step size, and $\mb{r}^{(i)}\in\mathbb{R}^{n}$ and $\delta\mb{x}^{(i)}\in\mathbb{R}^{n}$ denotes the residual form and the perturbation from the steady-state solution $\mb{x}_{\mathrm{tgt}}$ at the $i^{\mathrm{th}}$ time-instance, respectively.
In addition, the initial condition for this unsteady problem is that 
\begin{equation}
\delta\mb{x}^{(0)} = \delta\mb{x}_0.
\end{equation}
Writing \cref{eq:unsteady_control_2} in a compact form,



\begin{equation}
\label{eq:unsteady_control}
\mb{r}_{\mathrm{cl}}^{n_t}(\delta\mb{x}^{n_t}, \boldsymbol{\theta}, \mb{P})=
\begin{bmatrix}
\vdots \\
\mb{r}_{\mathrm{cl}}^{(i)} \\
\vdots
\end{bmatrix} = 0, \quad 
\delta\mb{x}^{n_t}
=
\begin{bmatrix}
\vdots \\
\delta\mb{x}^{(i)}  \\
\vdots
\end{bmatrix},
\end{equation}
where $i=1, \ldots , n_t$, and $\mb{r}_{\mathrm{cl}}^{n_t}$ and $\delta\mb{x}^{n_t}$ denotes the stacked residual and state variables, respectively.

\subsection{Analysis recap}

The function of interest consists of the sum of control and state costs,
\begin{equation}
\label{eq:cost_func}
\begin{aligned}
f &= \sum_{i=1}^{n_t} \left({\delta\mb{x}^{(i)}}^\intercal \mb{Q} \delta\mb{x}^{(i)} + {\delta\mb{x}^{(i)}}^\intercal\left(- \mb{P}^\intercal\mb{G}_{\mathrm{tgt}}\mb{S}^{-\intercal}  \right)\mb{S}\left(- \mb{S}^{-1} \mb{G}_{\mathrm{tgt}}^\intercal \mb{P}\right) \delta\mb{x}^{(i)}\right)\Delta t\\ 
&= \sum_{i=1}^{n_t} {\delta\mb{x}^{(i)}}^\intercal \left(\mb{Q} + \mb{P}^\intercal\mb{G}_{\mathrm{tgt}}\mb{S}^{-\intercal}\mb{G}_{\mathrm{tgt}}^\intercal \mb{P} \right) \delta\mb{x}^{(i)}\Delta t,
\end{aligned}
\end{equation}
and the residual form is defined by
\begin{equation}
\label{eq:res}
\mb{R}(\boldsymbol{\theta}, \mb{P}, \delta\mb{x}^{n_t}, \mb{d}) = 
\begin{bmatrix}
\check{\mb{r}}_{\mathrm{nl}}(\boldsymbol{\theta}, \mb{d})\\
\mb{R}_{\mathrm{ARE}}(\mb{P}, \mb{J}(\boldsymbol{\theta}, \mb{d}), \mb{G}(\boldsymbol{\theta}, \mb{d}))\\
\mb{r}_{\mathrm{cl}}^{n_t}(\delta\mb{x}^{n_t}, \boldsymbol{\theta}, \mb{P}, \mb{d})
\end{bmatrix}
=\mb{0}.
\end{equation}
The residual form $\mb{R}$ is not a matrix or a vector; it is a list containing all the individual residual forms written in either matrix or vector form.

\begin{figure*}[!ht]
\centering
\includegraphics[width=0.9\textwidth]{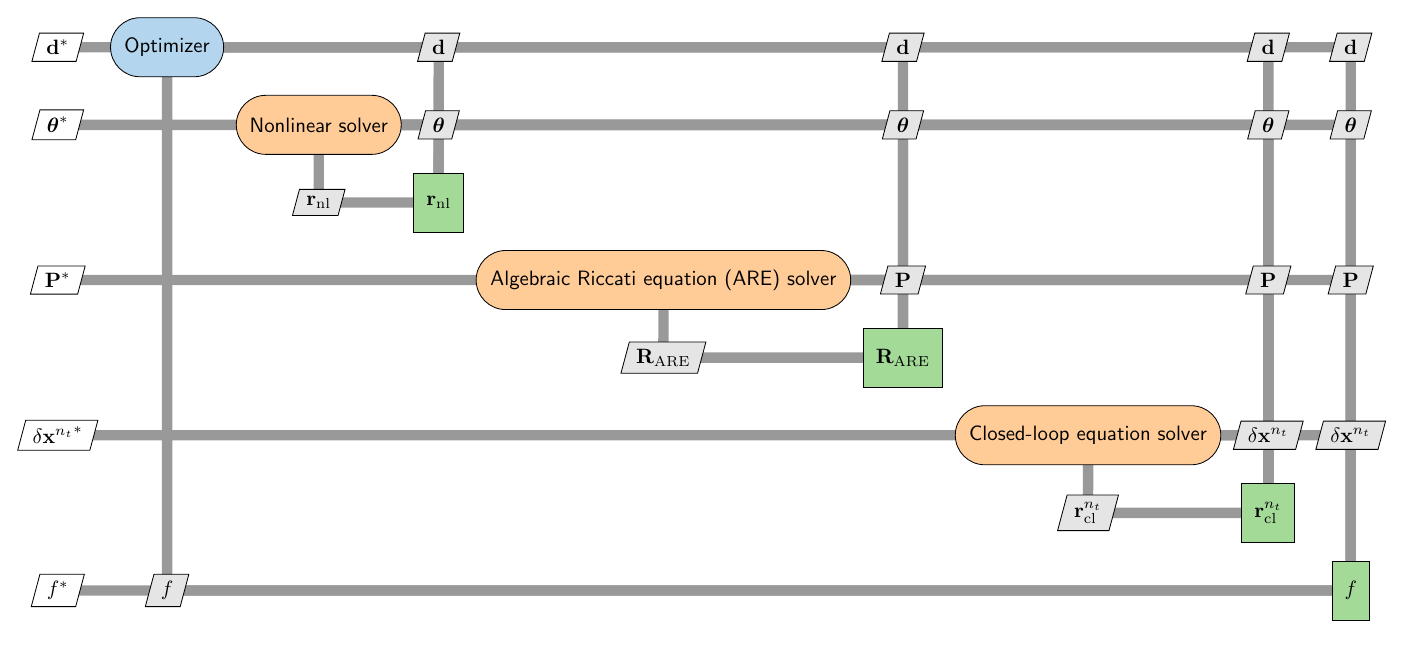}
\caption{XDSM for CCD analysis and optimization.}
\label{fig:XDSM}
\end{figure*}

The procedures discussed in this section are summarized in \cref{fig:XDSM} using an extended design structure matrix (XDSM) format~\cite{Lambe2012a}.
The figure shows that the information flows as follows.
First, the nonlinear solver computes the target state $\boldsymbol{\theta}$ by solving the steady-state residual $\mb{r}_{\mathrm{nl}}$.
Then, the ARE solver computes the Riccati solution $\mb{P}$ from the ARE residual $\mb{R}_{\mathrm{ARE}}$, which depends on $\boldsymbol{\theta}$ and the design variables $\mb{d}$.
Finally, the closed-loop equation solver integrates the unsteady residual $\mb{r}^{n_t}_{\mathrm{cl}}$ using the feedback law derived from $\mb{P}$.
The objective function $f$ is then evaluated from the resulting state trajectory.
The CCD formulation is summarized in \cref{alg:analysis}.

\begin{algorithm*}
    \caption{CCD Analysis}
    \begin{algorithmic}[1]
      \Procedure{CCD}{$\mb{d}, \hat{\mb{x}}_{\mathrm{tgt}}, \hat{\mb{u}}_{\mathrm{tgt}}, \mb{S}, \mb{Q}$} \Comment{Apply infinite horizon LQR to general nonlinear system.}
        \State $\boldsymbol{\theta} \rightarrow \check{\mb{r}}_{\mathrm{nl}}(\boldsymbol{\theta}; \hat{\mb{x}}_{\mathrm{tgt}}, \hat{\mb{u}}_{\mathrm{tgt}}, \mb{d}) = \mb{0}.$\Comment{Compute unknown target state and control using \cref{eq:steady_state_full_compact}.}
        \State $\mb{J}_{\mathrm{tgt}}=\f{\p\mb{r}_{\mathrm{nl}}}{\p \mb{x}}\left({\mb{x}_{\mathrm{tgt}}, \mb{u}_{\mathrm{tgt}}, \mb{d}}\right), \mb{G}_{\mathrm{tgt}}=\f{\p\mb{r}_{\mathrm{nl}}}{\p \mb{u}}\left({\mb{x}_{\mathrm{tgt}}, \mb{u}_{\mathrm{tgt}}, \mb{d}}\right)$\Comment{Linearize system around target state, as in \cref{eq:rlin}.}
        \State $\mb{P} \rightarrow \mb{R}_{\mathrm{ARE}}(\mb{P}; \mb{J}_{\mathrm{tgt}}, \mb{G}_{\mathrm{tgt}}, \mb{Q}, \mb{S}) = \mb{O}$\Comment{Solve ARE residual \cref{eq:res_ARE}.}
        \State $\mb{W} = - \mb{S}^{-1} \mb{G}_{\mathrm{tgt}}^\intercal \mb{P}$\Comment{Compute optimal state feedback matrix from \cref{eq:fsfbW}.}
        \State $\mb{r}_{\mathrm{cl}}^{n_t}(\delta\mb{x}^{n_t}, \boldsymbol{\theta}, \mb{P}, \mb{d}) = 0$\Comment{Integrate the closed-loop system to solve for the optimal trajectory in \cref{eq:unsteady_control_2}.}
        \State $f = f(\boldsymbol{\theta}, \mb{P}, \delta\mb{x}^{n_t}, \mb{d})$\Comment{Compute the total LQR cost defined in \cref{eq:cost_func}.}
        \State \textbf{return} $f$\Comment{Return the LQR cost over the optimal trajectory.}
      \EndProcedure
    \end{algorithmic}
    \label{alg:analysis}
  \end{algorithm*}

\section{Differentiable CCD}
\label{sec:oc-derivatives}

We derive the adjoint equation by forming the Lagrangian of the optimization problem, as detailed in \cref{ap:adjoint}.
The Lagrangian function is 
\begin{equation}
l = f + \boldsymbol{\psi}_{\mathrm{nl}}^\intercal\mb{r}_{\mathrm{nl}} + \mathrm{tr}\left(\boldsymbol{\Psi}_{\mathrm{ARE}}^\intercal\mb{R}_{\mathrm{ARE}}\right) + {\boldsymbol{\psi}_{\mathrm{cl}}^{n_t}}^\intercal\mb{r}_{\mathrm{cl}}^{n_t}.
\end{equation}
We have that
\begin{equation}
\label{eq:cadjoint_GE}
\begin{aligned}
\f{\p l}{\p \boldsymbol{\theta}} &= \f{\p f}{\p \boldsymbol{\theta}} + \boldsymbol{\psi}_{\mathrm{nl}}^\intercal \f{\p \mb{r}_{\mathrm{nl}}}{\p \boldsymbol{\theta}} + \f{\p\mathrm{tr}\left(\boldsymbol{\Psi}_{\mathrm{ARE}}^\intercal\mb{R}_{\mathrm{ARE}}\right)}{\p \boldsymbol{\theta}} + {\boldsymbol{\psi}_{\mathrm{cl}}^{n_t}}^\intercal\f{\p\mb{r}_{\mathrm{cl}}^{n_t}}{\p \boldsymbol{\theta}}=\mb{0},\\
\f{\p l}{\p \mb{P}} &= \f{\p f}{\p \mb{P}} + \underbrace{\boldsymbol{\psi}_{\mathrm{nl}}^\intercal \f{\p \mb{r}_{\mathrm{nl}}}{\p \mb{P}}}_{=\mb{O}} + \f{\p\mathrm{tr}\left(\boldsymbol{\Psi}_{\mathrm{ARE}}^\intercal\mb{R}_{\mathrm{ARE}}\right)}{\p \mb{P}} + {\boldsymbol{\psi}_{\mathrm{cl}}^{n_t}}^\intercal\f{\p\mb{r}_{\mathrm{cl}}^{n_t}}{\p \mb{P}}=\mb{O},\\
\f{\p l}{\p \delta\mb{x}^{n_t}} &= \f{\p f}{\p \delta\mb{x}^{n_t}} + \underbrace{\boldsymbol{\psi}_{\mathrm{nl}}^\intercal \f{\p \mb{r}_{\mathrm{nl}}}{\p \delta\mb{x}^{n_t}}}_{=\mb{0}} + \underbrace{\f{\p\mathrm{tr}\left(\boldsymbol{\Psi}_{\mathrm{ARE}}^\intercal\mb{R}_{\mathrm{ARE}}\right)}{\p \delta\mb{x}^{n_t}}}_{=\mb{0}} + {\boldsymbol{\psi}_{\mathrm{cl}}^{n_t}}^\intercal\f{\p \mb{r}_{\mathrm{cl}}^{n_t}}{\p \delta\mb{x}^{n_t}} = \mb{0} .\\
\end{aligned}
\end{equation}
These are the adjoint equations for the unsteady equation, ARE, and the steady-state equation, respectively.  

The partial derivatives for the higher level residual only depend on the state variables from the same level or the lower levels; for example, $\mb{R}_{\mathrm{ARE}}$ is dependent on $\mb{x}_{\mathrm{tgt}}$ and $\mb{P}$, but it is independent of $\delta\mb{x}^{n_t}$.
Leveraging this special structure, we solve the coupled adjoint equation separately by solving each adjoint equation, starting from the top level and working downwards, using block back-substitution as follows: 
\begin{equation}
\begin{aligned}
{\boldsymbol{\psi}_{\mathrm{cl}}^{n_t}} &\rightarrow {\boldsymbol{\psi}_{\mathrm{cl}}^{n_t}}^\intercal\f{\p \mb{r}_{\mathrm{cl}}^{n_t}}{\p \delta\mb{x}^{n_t}} = - \f{\p f}{\p \delta\mb{x}^{n_t}},\\
\boldsymbol{\Psi}_{\mathrm{ARE}} &\rightarrow \f{\p\mathrm{tr}\left(\boldsymbol{\Psi}_{\mathrm{ARE}}^\intercal\mb{R}_{\mathrm{ARE}}\right)}{\p \mb{P}} = - \f{\p l}{\p \mb{P}} - {\boldsymbol{\psi}_{\mathrm{cl}}^{n_t}}^\intercal\f{\p\mb{r}_{\mathrm{cl}}^{n_t}}{\p \mb{P}},\notag\\
\boldsymbol{\psi}_{\mathrm{nl}}^\intercal \f{\p \mb{r}_{\mathrm{nl}}}{\p \boldsymbol{\theta}} &\rightarrow \boldsymbol{\psi}_{\mathrm{nl}}^\intercal \f{\p \mb{r}_{\mathrm{nl}}}{\p \boldsymbol{\theta}} = -\f{\p f}{\p \boldsymbol{\theta}} - \f{\p\mathrm{tr}\left(\boldsymbol{\Psi}_{\mathrm{ARE}}^\intercal\mb{R}_{\mathrm{ARE}}\right)}{\p \boldsymbol{\theta}} - {\boldsymbol{\psi}_{\mathrm{cl}}^{n_t}}^\intercal\f{\p\mb{r}_{\mathrm{cl}}^{n_t}}{\p \boldsymbol{\theta}}.\notag
\end{aligned}
\end{equation}
By avoiding the direct solution of the coupled adjoint equation, we reduce computational cost by using the efficient adjoint solver for each module. 

After the adjoint variables are computed, we can use the following formula to evaluate the total derivative
\begin{equation}
\label{eq:total_derivative_intro}
\resizebox{.8\hsize}{!}{$\f{\d f}{\d \mb{d}} = \f{\p f}{\p \mb{d}} + \boldsymbol{\psi}_{\mathrm{nl}}^\intercal \f{\p \mb{r}_{\mathrm{nl}}}{\p \mb{d}} + \f{\p  \mathrm{tr}\left(\boldsymbol{\Psi}_{\mathrm{ARE}}^\intercal\mb{R}_{\mathrm{ARE}}\right)}{\p \mb{d}} + {\boldsymbol{\psi}_{\mathrm{cl}}^{n_t}}^\intercal\f{\p \mb{r}_{\mathrm{cl}}^{n_t}}{\p \mb{d}}.$}
\end{equation}
We detail the solution of the three adjoint equations in this section.
In the following discussion, rather than a special right-hand side (RHS), we consider a general RHS.
Then, we discuss the analytic formulas for partial derivatives for the RHS.

\subsection{Closed-loop system adjoint equation solution}

The unsteady adjoint equation is written as
\begin{equation}
 \f{\p \mb{r}_{\mathrm{cl}}^{n_t}}{\p \delta\mb{x}^{n_t}}^\intercal{\boldsymbol{\psi}_{\mathrm{cl}}^{n_t}} = - \f{\p f}{\p \delta\mb{x}^{n_t}}^\intercal.
\label{eq:closed_loop_adjoint}
\end{equation}
Solving this linear equation directly is inefficient due to its large size.
Instead, we solve step by step in the reverse direction.
In \cref{sec:unsteady_adjoint}, we show an example with a forward Euler scheme.
Other schemes can be implemented in a similar way.
Wang et al.~\cite{Wang2009} and Zhang et al.~\cite{Zhang2022a} provide more details about efficiently solving unsteady adjoint equations.

\subsection{Adjoint Algebraic Riccati equation (ARE) solution}

The ARE can be differentiated efficiently using the adjoint method.
Following the derivation by Kao and Hennequin~\cite{Kao2020a}.
We solve the following matrix adjoint equation, which is categorized as a Lyapunov equation defined by
\begin{equation}
\tilde{\mb{J}}\boldsymbol{\Psi}_{\mathrm{ARE}} + \boldsymbol{\Psi}_{\mathrm{ARE}}\tilde{\mb{J}}^\intercal + \f{1}{2}\left(\left(\f{\p f}{\p \mb{P}} + \f{\p \mb{r}_{\mathrm{cl}}^{n_t}}{\p \mb{P}}^\intercal\boldsymbol{\psi}_{\mathrm{cl}}^{n_t}\right) + \left(\f{\p f}{\p \mb{P}} + \f{\p \mb{r}_{\mathrm{cl}}^{n_t}}{\p \mb{P}}^\intercal\boldsymbol{\psi}_{\mathrm{cl}}^{n_t}\right)^\intercal\right) = \mb{O},
\label{eq:ARE_adjoint_equation}
\end{equation}
where $\boldsymbol{\Psi}_{\mathrm{ARE}}$ is the adjoint matrix to be computed, $\tilde{\mb{J}}$ is defined by
\begin{equation}
\tilde{\mb{J}} = \mb{J}_{\mathrm{tgt}} - \mb{G}_{\mathrm{tgt}}\mb{S}^{-1}\mb{G}_{\mathrm{tgt}}^\intercal\mb{P}.
\end{equation}
In \cref{sec:pr_pP_1}, we give an example for ${\p \mb{r}_{\mathrm{cl}}^{n_t}}/{\p \mb{P}}^\intercal\boldsymbol{\psi}_{\mathrm{cl}}^{n_t}$.

\subsection{Adjoint nonlinear steady-state equation solution}

The adjoint equation of the steady-state equation is defined by
\begin{equation}
\resizebox{.85\hsize}{!}{$\f{\p \check{\mb{r}}_{\mathrm{nl}}}{\p \boldsymbol{\theta}}^\intercal \boldsymbol{\psi}_{\mathrm{nl}} = -\left(\f{\p f}{\p \boldsymbol{\theta}} + \boldsymbol{\psi}_{\mathrm{cl}}^\intercal\f{\p\mb{r}_{\mathrm{cl}}^{n_t}}{\p \boldsymbol{\theta}} + \f{\p\mathrm{tr}\left(\boldsymbol{\Psi}_{\mathrm{ARE}}^\intercal\mb{R}_{\mathrm{ARE}}\right)}{\p \boldsymbol{\theta}}\right).$}
\label{eq:nonlin_steady_state_adjoint}
\end{equation}
The matrices and vectors shown in the adjoint equation can be constructed directly or using the partial derivatives with respect to $\mb{x}_{\mathrm{tgt}}$ and $\mb{u}_{\mathrm{tgt}}$.
${\p \check{\mb{r}}_{\mathrm{nl}}}/{\p \boldsymbol{\theta}}$ can be contructed using $\mb{J}_{\mathrm{tgt}}$ and $\mb{G}_{\mathrm{tgt}}$.
${\p f}/{\p \boldsymbol{\theta}}$ can be constructed using ${\p f}/{\p \mb{x}_{\mathrm{tgt}}}$ and ${\p f}/{\p \mb{u}_{\mathrm{tgt}}}$.
In \cref{sec:pr_pP_2}, we give an example on how to construct $\boldsymbol{\psi}_{\mathrm{cl}}^\intercal {\p\mb{r}_{\mathrm{cl}}^{n_t}}/{\p \boldsymbol{\theta}}$ and ${\p\mathrm{tr}\left(\boldsymbol{\Psi}_{\mathrm{ARE}}^\intercal\mb{R}_{\mathrm{ARE}}\right)}/{\p \boldsymbol{\theta}}$.

\subsection{Total derivative}

The total derivative equation was presented earlier in \cref{eq:total_derivative_intro}; we rewrite it here and expand its terms. 
\begin{equation}
\f{\d f}{\d \mb{d}} = \f{\p f}{\p \mb{d}} + \boldsymbol{\psi}_{\mathrm{nl}}^\intercal \f{\p \mb{r}_{\mathrm{nl}}}{\p \mb{d}} + \f{\p  \mathrm{tr}\left(\boldsymbol{\Psi}_{\mathrm{ARE}}^\intercal\mb{R}_{\mathrm{ARE}}\right)}{\p \mb{d}} + {\boldsymbol{\psi}_{\mathrm{cl}}^{n_t}}^\intercal\f{\p \mb{r}_{\mathrm{cl}}^{n_t}}{\p \mb{d}}.
\label{eq:total_derivative}
\end{equation}
The term ${\p  \mathrm{tr}\left(\boldsymbol{\Psi}_{\mathrm{ARE}}^\intercal\mb{R}_{\mathrm{ARE}}\right)}/{\p \mb{d}}$ can be further expanded as
\begin{equation}
\label{eq:tr_der}
\f{\p\mathrm{tr}\left(\boldsymbol{\Psi}_{\mathrm{ARE}}^\intercal\mb{R}_{\mathrm{ARE}}\right)}{\p \mb{d}} = 
\f{\p \mathrm{tr}\left(\mb{J}_{\mathrm{tgt}}^\intercal\overbar{\mb{J}}_{\mathrm{tgt}}\right)}{\p \mb{d}} + \f{\p \mathrm{tr}\left(\mb{G}_{\mathrm{tgt}}^\intercal \overbar{\mb{G}}_{\mathrm{tgt}}\right)}{\p \mb{d}} 
\end{equation}
where the seeds $\overbar{\mb{J}}_{\mathrm{tgt}}$ and $\overbar{\mb{G}}_{\mathrm{tgt}}$ are defined by
\begin{equation}
\begin{aligned}
\overbar{\mb{J}}_{\mathrm{tgt}} &= \mb{P}\boldsymbol{\Psi}_{\mathrm{ARE}}^\intercal + \mb{P}^\intercal \boldsymbol{\Psi}_{\mathrm{ARE}},\\
\overbar{\mb{G}}_{\mathrm{tgt}} &= -\mb{P}^\intercal \boldsymbol{\Psi}_{\mathrm{ARE}} \mb{P}^\intercal\mb{G}_{\mathrm{tgt}}\mb{S}^{-\intercal}-\mb{P} \boldsymbol{\Psi}_{\mathrm{ARE}}^\intercal \mb{P}\mb{G}_{\mathrm{tgt}}\mb{S}^{-1}.
\end{aligned}
\end{equation}
The matrices $\overbar{\mb{J}}_{\mathrm{tgt}}, \overbar{\mb{G}}_{\mathrm{tgt}}$ are treated as constant weighting factors when evaluating derivatives with respect to $\mb{d}$ in \cref{eq:tr_der}. 
For more details about ${\boldsymbol{\psi}_{\mathrm{cl}}^{n_t}}^\intercal \left({\p \mb{r}_{\mathrm{cl}}^{n_t}}/{\p \mb{d}}\right)$, see \cref{sec:pr_pP_3}.

The procedure for calculating the total derivatives is shown in \cref{alg:derivatives}.

\begin{algorithm}
\caption{CCD Derivative Computation}
\begin{algorithmic}[1]
\Procedure{$\frac{\d f}{\d \mb{d}}$}{$\mb{d}, \boldsymbol{\theta}, \mb{{P}}, \delta\mb{x}^{n_t}$} \Comment{Compute derivative of function of interest w.r.t. design variables.}
\State $\boldsymbol{\psi}_{\mathrm{cl}}^{n_t} \rightarrow {\boldsymbol{\psi}_{\mathrm{cl}}^{n_t}}^\intercal\f{\p \mb{r}_{\mathrm{cl}}^{n_t}}{\p \delta\mb{x}^{n_t}}(\boldsymbol{\theta}, \mb{P}, \delta \mb{x}^{n_t}, \mb{d}) = - \f{\p f}{\p \delta\mb{x}^{n_t}}(\boldsymbol{\theta}, \mb{P}, \delta \mb{x}^{n_t}, \mb{d})$\Comment{Solve for the unsteady adjoint for closed-loop system using \cref{eq:closed_loop_adjoint}.}
\State $\boldsymbol{\Psi}_{\mathrm{ARE}} \rightarrow \f{\p\mathrm{tr}\left(\boldsymbol{\Psi}_{\mathrm{ARE}}^\intercal\mb{R}_{\mathrm{ARE}}\right)}{\p \mb{P}}(\boldsymbol{\theta}, \mb{P}, \mb{d}) = - \f{\p l}{\p \mb{P}} - {\boldsymbol{\psi}_{\mathrm{cl}}^{n_t}}^\intercal\f{\p\mb{r}_{\mathrm{cl}}^{n_t}}{\p \mb{P}}$ \Comment{Solve for the ARE adjoint via \cref{eq:ARE_adjoint_equation}.}
\State $\boldsymbol{\psi}_{\mathrm{nl}}  \rightarrow \boldsymbol{\psi}_{\mathrm{nl}}^\intercal \f{\p \mb{r}_{\mathrm{nl}}}{\p \boldsymbol{\theta}}(\boldsymbol{\theta},\mb{d}) = -\f{\p f}{\p \boldsymbol{\theta}} - \f{\p\mathrm{tr}\left(\boldsymbol{\Psi}_{\mathrm{ARE}}^\intercal\mb{R}_{\mathrm{ARE}}\right)}{\p \boldsymbol{\theta}} - {\boldsymbol{\psi}_{\mathrm{cl}}^{n_t}}^\intercal\f{\p\mb{r}_{\mathrm{cl}}^{n_t}}{\p \boldsymbol{\theta}}$\Comment{Solve for the steady-state adjoint equation by \cref{eq:nonlin_steady_state_adjoint}.}
\State $\f{\d f}{\d \mb{d}} = \f{\p f}{\p \mb{d}} + \boldsymbol{\psi}_{\mathrm{nl}}^\intercal \f{\p \mb{r}_{\mathrm{nl}}}{\p \mb{d}} + \f{\p  \mathrm{tr}\left(\boldsymbol{\Psi}_{\mathrm{ARE}}^\intercal\mb{R}_{\mathrm{ARE}}\right)}{\p \mb{d}} + {\boldsymbol{\psi}_{\mathrm{cl}}^{n_t}}^\intercal\f{\p \mb{r}_{\mathrm{cl}}^{n_t}}{\p \mb{d}}$\Comment{Compute total derivative using \cref{eq:total_derivative}.}
\State \textbf{return} $\f{\d f}{\d \mb{d}}$
\EndProcedure
\end{algorithmic}
\label{alg:derivatives}
\end{algorithm}

\subsection{Function of interest}

To compute the total derivatives, we also need to compute several partial derivatives of the function of interest, including $\p f / \ \mb{x}, \p f / \p \delta\mb{x}^{(i)}, \p f / \p \mb{P}, \p f / \p \boldsymbol{\theta}$ and $\p f / \p \mb{G}_{\mathrm{tgt}}$.
They are straightforward to evaluate. 
In \cref{sec:FoI_partials}, we present results for the case of the cost function defined by \cref{eq:cost_func}.


\section{Numerical example}
\label{sec:numerical_example}

\subsection{Derivative computation verification}
\subsubsection{Problem setup}
\begin{figure}[htbp]
\centering
\includegraphics[width=0.9\columnwidth]{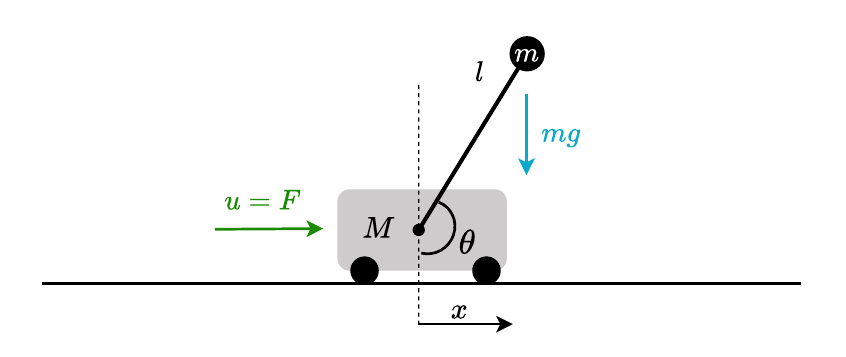}
\caption{The control input $u$ in the cart-pole problem causes the cart to slide with no friction.}
\label{fig:cart}
\end{figure}

We consider the classic inverted cart-pole problem shown in figure \cref{fig:cart} with the parameters defined by \cref{tab:cart}.
The equilibrium point in this problem is directly specified and therefore the nonlinear steady-state equation \cref{eq:steady_state_full_compact} 
does not need to be solved. 
Thus, the problem is only composed of the optimal control module and the closed-loop system.

Following \cref{eq:dyn_open}, the dynamical equation for the inverted pole on a cart is defined by
\begin{equation}
\mb{x} = \mb{r}_{\mathrm{nl}} + \mb{G}\mb{u} ,
\end{equation}
where
\begin{equation}
\begin{gathered}
\mb{x} = 
\begin{bmatrix}
x \\
v \\
\theta \\
\omega
\end{bmatrix},\quad
\mb{r}_{\mathrm{nl}} =
\begin{bmatrix}
v \\
\f{-m^2L^2g\cos{\theta}\sin{\theta} + mL^2(mL\omega^2\sin{\theta} - \delta v)}{mL^2(M + m(1 - \cos^2{\theta}))} \\
\omega \\
\f{(m + M)mgL\sin{\theta} - mL\cos{\theta}(mL\omega^2\sin{\theta} - \delta v)}{mL^2(M + m(1 - \cos^2{\theta}))}
\end{bmatrix},\\
\mb{G}=
\f{1}{L(M + m(1 - \cos^2{\theta}))}
\begin{bmatrix}
0 \\
L  \\
0 \\
\cos{\theta} 
\end{bmatrix}.
\end{gathered}
\end{equation}
These equations are written in dimensionless form, and the parameters are defined by \cref{tab:cart}.
The following two fixed-point solutions exist for this problem: 
\begin{equation}
\mb{x}_{\mathrm{tgt}} = 
\begin{bmatrix}
x^* \\
0 \\
0 \\
0
\end{bmatrix}, \quad \mathrm{and}\quad
\mb{x}_{\mathrm{tgt}} = 
\begin{bmatrix}
x^* \\
0 \\
\pi \\
0
\end{bmatrix},\quad
\mb{u}_{\mathrm{tgt}} = 0,
\label{eq:stable_points} 
\end{equation}
where the $x$ coordinates $x^*$ is arbitrary.
The target is to stabilize the system with the pole pointing up ($\theta = \pi$), and we pick $x^* = 1$ to avoid ambiguity.

\begin{table}[ht]
\centering
\caption{Parameters for the cart-pole problem.}
\begin{tabular}{lr}
\toprule
Parameter & Value \\
\midrule
$m$ & 1 \\
$M$ & 5 \\
$L$ & 2 \\
$g$ & 10 \\
\bottomrule
\end{tabular}
\label{tab:cart}
\end{table}

\subsubsection{Total derivatives}
The total derivatives of the LQR cost function computed using \cref{alg:derivatives} and are compared to finite difference results in \cref{tab:der_verification}.
We obtained a 5--7 digit match between the adjoint method and the finite differences (using forward differences with a step size of $\times 10^{-6}$).
This verifies our adjoint method implementation.
\begin{table}[ht]
\centering
\small\setlength\tabcolsep{4.5pt}
\caption{Verification of the derivatives $\d f / \d \mb{d}$.}
\begin{tabular}{lccc}
\toprule
& $m$ & $M$ & $L$ \\
\midrule
Adjoint & $716.01\underline{2631}$ & $766.693\underline{430}$ & $172.7323\underline{23}$ \\
FD & $716.01\underline{3365}$ & $766.693\underline{265}$ & $172.7323\underline{30}$ \\
\bottomrule
\end{tabular}
\label{tab:der_verification}
\end{table}

\subsection{Cart-pole problem}
\subsubsection{Problem setup}
We now consider a CCD problem for the cart-pole described in the previous section, where the goal is to optimize the parameters ($\mb{d} = \begin{bmatrix}m, M, L\end{bmatrix}^\intercal$) of the system to balance the pendulum at the equilibrium point $\mb{x}_{\mathrm{tgt}}$ with the lowest LQR cost. 
The baseline case for this example uses the parameters shown in \cref{tab:cart}.
The initial state $\mb{x}_{\mathrm{init}}$ is a randomly chosen state close enough to $\mb{x}_{\mathrm{tgt}}$ for LQR satisfaction and is given as
\begin{equation}
\mb{x}_{\mathrm{tgt}} = 
\begin{bmatrix}
0 \\
0 \\
\pi \\
0
\end{bmatrix},
\quad
\mb{x}_{\mathrm{init}} = 
\begin{bmatrix}
-1 \\
0 \\
2 \\
0
\end{bmatrix}
.
\end{equation}
We added a lower bound on the total mass and bounds on the mass of the cart and the pole, as well as bounds on the length of the pole. 
Bounds on the design variables are added, including a minimum mass constraint to prevent trivially small designs, 
\begin{equation}
\begin{bmatrix}
0.5 \\
2.5 \\
1 
\end{bmatrix} \quad \leq\quad
\begin{bmatrix}
m \\
M \\
L 
\end{bmatrix} \quad \leq\quad
\begin{bmatrix}
2.0 \\
7.5 \\
2
\end{bmatrix}
, \quad m + M \geq 3.5.
\end{equation}
The baseline case shown in \cref{tab:cart} satisfies the above bounds. 

\subsubsection{Optimization results}
\begin{figure}[htbp]
\centering
\includegraphics[width=0.9\columnwidth]{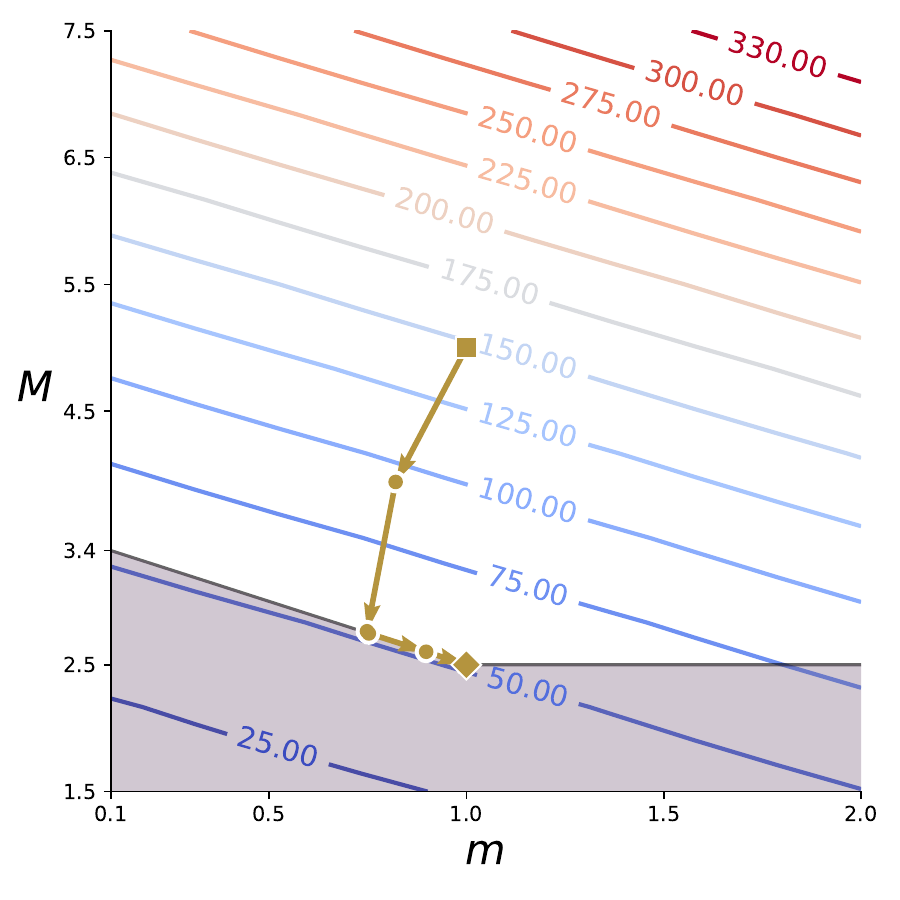}
\caption{Cost contour of $M$ and $m$ with $l = 1$, where the infeasible region is shaded.}
\label{fig:cart_contour}
\end{figure}

The LQR optimization problem for this example is as follows:
\begin{equation}
\label{eq:cartpole-opt-problem}
\begin{aligned}
\min_{\mb{d}} & \quad \sum_{i=1}^{n_t} {\delta\mb{x}^{(i)}}^\intercal \left(\mb{Q} + \mb{P}^\intercal\mb{G}_{\mathrm{tgt}}\mb{S}^{-\intercal}\mb{G}_{\mathrm{tgt}}^\intercal \mb{P} \right) \delta\mb{x}^{(i)}\Delta t, \\
\text{subject to} & \quad \mb{d}_{\mathrm{min}} \leq \mb{d} \leq \mb{d}_{\mathrm{max}} , \\
& \quad m + M \geq 3.5, \\
& \quad \mb{R}(\boldsymbol{\theta}, \mb{P}, \delta\mb{x}^{n_t}, \mb{d})=\mb{0}, \\
\end{aligned}   
\end{equation}
where $\mb{d}_{\mathrm{min}}$ and $\mb{d}_{\mathrm{max}}$ are the upper and lower bounds on the design variables, and $\mb{R}(\boldsymbol{\theta}, \mb{P}, \delta\mb{x}^{n_t}, \mb{d})$ is the residual form defined by \cref{eq:res}.
We also set $\mb{Q} =0.1\mb{I}$ and $\mb{S = I}$.
Arbitrary positive semidefinte matrices $\mb{P}$ and $\mb{R}$ can be used in the cost function.
These values can be used to reflect the importance of control and state error.
This optimization problem is solved using the adjoint method framework to calculate the cost sensitivity of the design variables. 
We use SNOPT~\cite{Gill2005a} through the pyOptSparse wrapper~\cite{Wu2020a} to solve the optimization problem starting with the baseline design variables. 
This results in the optimal design shown in \cref{tab:cart_opt}. 
A contour of the cost function with respect to $M$ and $m$ is shown in \cref{fig:cart_contour} with a fixed pole length of $l = 1$. 
The optimization path can be seen, as well as the infeasible region resulting from the design constraints.

\begin{table}[ht]
\centering
\caption{Optimal cart-pole parameters}
\begin{tabular}{lcr}
\toprule
Parameter & Baseline & Optimal Value \\
\midrule
$m$ & 1.0 & 1.0\\
$M$ & 5.0 & 2.5\\
$L$ & 2.0 & 1.0\\
\bottomrule
\end{tabular}
\label{tab:cart_opt}
\end{table}

\begin{figure}[htbp]
\centering
\includegraphics[width=0.9\columnwidth]{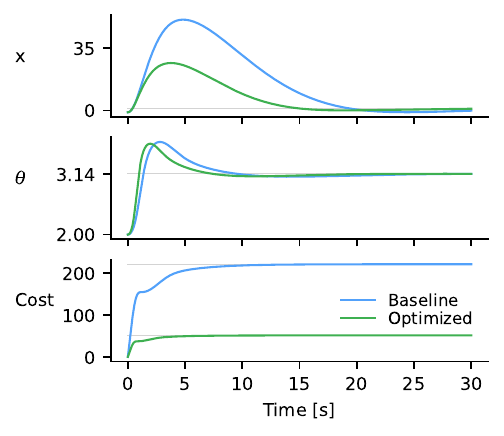}
\caption{Comparison of optimized and baseline results.}
\label{fig:optimization}
\end{figure}

The time response of the optimized system is compared to the baseline case in \cref{fig:optimization}, which shows a $76.79\%$ LQR cost reduction. 
In the case of the cart-pole, the optimization result is expected because the cost of control decreases as the values of the design variables decrease (hence the need for bounds on the design variables). 
However, in a more complex optimization problem such as that discussed in~\cref{subsect:quad}, our algorithm can find non-obvious solutions. 


\subsection{Quadrotor problem}
\label{subsect:quad}
\subsubsection{Problem setup}
In the second example, we solve a CCD optimization problem of a quadrotor drone, where we concurrently optimize the drone's LQR controller and the rotor blade geometry.
The LQR controller aims to bring the quadrotor state to steady hover starting from a perturbed state.

\begin{figure}[htbp]
\centering
\includegraphics[width=0.7\columnwidth]{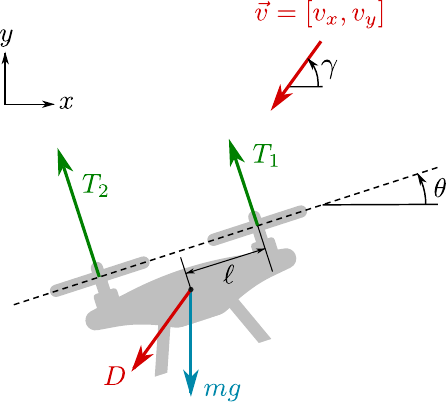}
\caption{Schematic of the quadrotor's longitudinal dynamics.}
\label{fig:quadrotor-schematics}
\end{figure}

We consider a three-degree-of-freedom longitudinal motion in the vertical plane as shown in \cref{fig:quadrotor-schematics}, where $x$, $y$, and $\theta$ are the drone's horizontal location, altitude, and pitch angle, respectively.
We also define the corresponding linear and angular velocities as $v_x$, $v_y$, and $\omega$.
We assumed that the yaw attitude is fixed at the configuration, and each of the right two rotors produces the thrust $T_1$, and each of the left two rotors produces $T_2$.

The longitudinal equations of motion are
\begin{equation}
\begin{aligned}
m \dot{v}_x &= -2(T_1 + T_2) \sin \theta - D \cos \gamma ~, \\
m \dot{v}_y &= 2(T_1 + T_2) \cos \theta - D \sin \gamma - mg ~, \\
I \dot{\omega} &= 2(T_1 - T_2) \ell ~,
\end{aligned}
\end{equation}
where $m$, $I$, and $\ell$ are the quadrotor mass, moment of inertia, and the distance between thrust acting points and the center of gravity, respectively.
The flight path angle $\gamma$ is related to the velocity by
\begin{equation}
\cos \gamma = \frac{v_x}{\sqrt{v_x^2 + v_y^2}}, \quad \sin \gamma = \frac{v_y}{\sqrt{v_x^2 + v_y^2}} ~.
\end{equation}
The body drag $D$ is given by
\begin{equation}
D = \frac{1}{2} \rho S C_D (v_x^2 + v_y^2) =: \beta (v_x^2 + v_y^2)~,
\end{equation}
where $\rho$ is the air density, $S$ is the body frontal area, and $C_D$ is the drag coefficient.
We introduced a new constant $\beta := 1/2 \rho S C_D$ to simplify the notation.
\Cref{table:quadrotor_parameters} summarizes the parameter values used in this example based on a notional drone presented in Quan et al.~\cite[Ch.~3.2]{Quan2020a}.

\begin{table}[h]
\centering
\caption{Parameters for the quadrotor problem.}
\begin{tabular}{l l l}
\toprule
Parameter & Value & Units\\
\midrule
$m$ & 1.4 & \unit{kg}\\
$I$ & 0.0211 & \unit{kg.m^2} \\
$\ell$ & 0.159 & \unit{m}\\
$\beta$ & 0.1365 & \unit{N / (m/s)^2} \\ 
\bottomrule
\end{tabular}
\label{table:quadrotor_parameters}
\end{table}

At every time instance, we run a blade element momentum (BEM) analysis to compute the rotor thrust $T$ and power $P$ as a function of the rotor speed $\Omega$ and rotor blade design $\mb{d}$.
The design variables consist of the spanwise distribution of the blade's chord and twist, parametrized by b-splines with 10 spanwise control points; therefore, the total number of plant design variables is 20.
The rotor diameter is fixed at 0.266~m~\cite{Quan2020a}.
We use CCBlade for BEM analysis~\cite{Ning2021} and wrap it using OpenMDAO~\cite{Gray2019a}. 
The derivatives of $T$ and $P$ with respect to $\Omega, \mb{d}$ are computed using OpenMDAO's coupled adjoint method~\cite{Hwang2018a} with partial derivatives from algorithmic differentiation~\cite{Martins2013a}. 
 The coupled derivatives are passed to the framework introduced in \cref{sec:oc-derivatives}.

The state variables are defined as $\mb{x} = [x,y,\theta,v_x,v_y,\omega]^T$ and the control variables are $\mb{u} = [\Omega_1,~\Omega_2]$.
The dynamical equations are
\begin{equation}
\mb{\dot{x}} =
\begin{bmatrix}
\dot{x} \\ \dot{y} \\ \dot{\theta} \\ \dot{v_x}  \\ \dot{v_y}  \\ \dot{\omega}
\end{bmatrix}
=
\begin{bmatrix}
v_x \\
v_y \\
\omega \\
\left( -2(T_1 + T_2) \sin \theta - \beta v_x \sqrt{v_x^2 + v_y^2} \right) / m \\
\left( 2(T_1 + T_2) \cos \theta - \beta v_y \sqrt{v_x^2 + v_y^2} \right) / m - g \\
2 (T_1 - T_2) \ell / I 
\end{bmatrix}
= \mb{r}_{\mathrm{nl}}(\mb{x}, \mb{u}, \mb{d}) ~.
\end{equation}
Because we aim to bring the quadrotor to a steady hovering state, $\theta=v_x=v_y=\omega=0$ at the target state.
The $x$ and $y$ values at the target can be any value, but we set them to zero for convenience.
This gives the target state vector of $\mb{x}_{\mathrm{tgt}} = \mb{0}$.
The control input at the target state, $\mb{u}_{\mathrm{tgt}}$, is unknown and dependent on the rotor design because the rotor speed at steady hover varies as the optimizer changes rotor blade design.
The target control input is determined by solving \cref{eq:steady_state_full} given $\mb{d}$ at every optimization iteration.
For the initial perturbed states, we use $\mb{x}^{(0)} = [1,1,0.1,0.5,0.3,0.05]$.

\subsubsection{Optimization Problem}

In CCD, optimizing the plant design for a dynamic performance metric, such as the LQR cost, often results in deteriorating a steady performance metric and vice versa.
This means that CCD often leads to a multiobjective optimization problem with two competing objectives: dynamic and steady performance metrics.

In this quadrotor example, we use the LQR's quadratic cost as a dynamic performance metric, and the power required for steady hover as a steady performance metric.
We solve this multiobjective optimization problem using the epsilon-constraint method~\cite[Ch.~9.3.2]{Martins2022}.
The optimization problem is as follows:
\begin{align}
\begin{aligned}
\min_{\mb{d}} & \quad \sum_{i=1}^{n_t} {\delta\mb{x}^{(i)}}^\intercal \left(\mb{Q} + \mb{P}^\intercal\mb{G}_{\mathrm{tgt}}\mb{S}^{-\intercal}\mb{G}_{\mathrm{tgt}}^\intercal \mb{P} \right) \delta\mb{x}^{(i)}\Delta t ~, \\
\text{subject to} & \quad P_{\text{hover}}(\mb{u}_{\mathrm{tgt}}, \mb{d}) \leq (1 + \epsilon) P_{\text{hover, min}} ~,\\
& \quad \mb{R}(\boldsymbol{\theta}, \mb{P}, \delta\mb{x}^{n_t}, \mb{d})=\mb{0},
\label{eq:quadrotor-opt-problem}
\end{aligned}   
\end{align}
where $P_{\mathrm{hover}}$ is the power required for steady hover.
This is a steady metric in the sense that it is independent of the state trajectory $\mb{x}(t)$ and control history $\mb{u}(t)$ for attitude recovery.
On the right-hand side of the constraint, $P_{\text{hover, min}}$ is the minimum value for the steady hovering power, and $\epsilon$ is a relaxation factor.
The objective is approximated by a numerical integration across a finite horizon following \cref{eq:cost_func}.
We set $\mb{Q}=\mb{I}$ and $\mb{S}=0.01\mb{I}$.
For the blade design variables, we bound the chord to be between 5\% and 35\% to avoid unreasonable chord lengths.
The twist variables are unbounded.  

Before solving the CCD problem in \cref{eq:quadrotor-opt-problem}, we solve the single-objective blade design optimization to minimize $P_{\mathrm{hover}}$ by varying the blade design without considering optimal control.
This yields the value of $P_{\text{hover, min}}$, which we use for the CCD optimizations.
Similarly to the cart-pole example, we use SNOPT~\cite{Gill2005a} and pyOptSparse~\cite{Wu2020a} to solve the optimization problems.

\subsubsection{Results}

\begin{figure}[htbp]
    \centering
    \includegraphics[width=1.0\columnwidth]{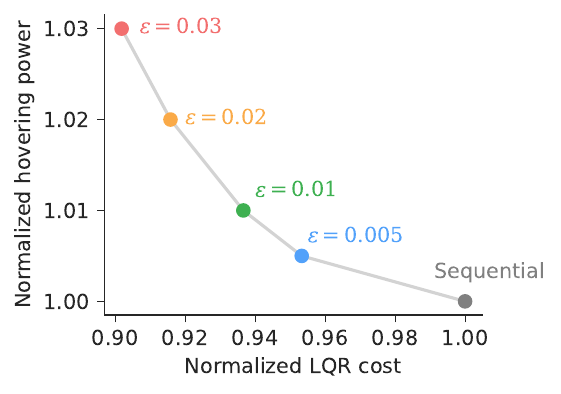}
    \caption{Pareto front representing the hovering power and LQR cost trade-off.}
    \label{fig:quadrotor-pareto-front}
\end{figure}

\Cref{fig:quadrotor-pareto-front} shows the Pareto front representing the trade-off between the steady hovering power and the LQR cost.
The right-most point shows the results of sequential optimization, where we first optimize the blade design for minimum hovering power and then solve the fixed-design optimal control to compute the LQR cost.
This design gives the minimum hovering power; however, its LQR cost is the worst because the blade design decision does not account for the attitude recovery control.
The other points on the Pareto front correspond to the CCD optimizations with $\epsilon=0.005, 0.01, 0.02, 0.03$.
For example, $\epsilon=0.03$ means that it minimizes the LQR cost while allowing 3\% higher power for steady hovering.

\begin{figure}[htbp]
    \centering
    \includegraphics[width=0.9\columnwidth]{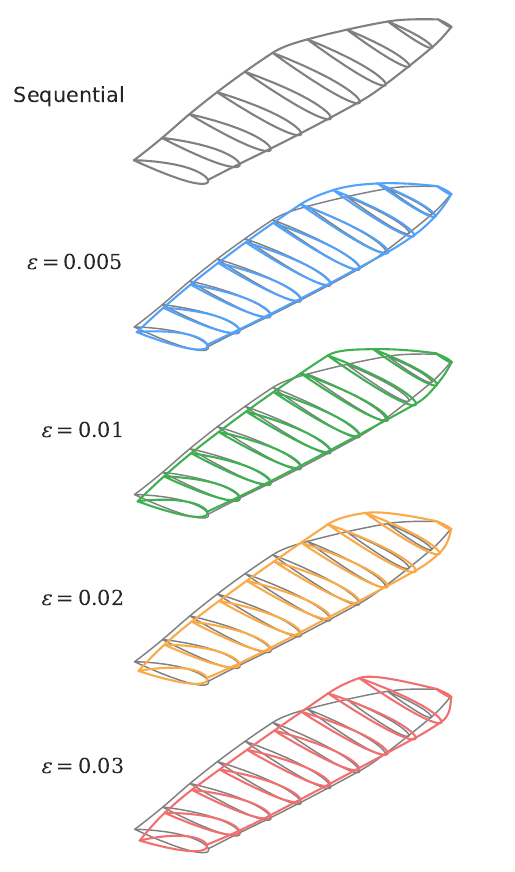}
    \caption{Optimal rotor blade geometries with sequential optimization and CCD with different $\epsilon$ values.}
    \label{fig:quadrotor-rotor-design-3D}
\end{figure}

\begin{figure}[htbp]
    \centering
    \includegraphics[width=1.0\columnwidth]{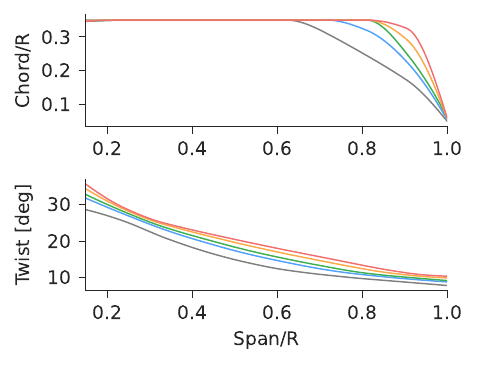}
    \caption{Comparison of the rotor's chord and twist distributions, where each line corresponds to a point in \cref{fig:quadrotor-pareto-front} with the same color.}
    \label{fig:quadrotor-rotor-design}
\end{figure}

As we prioritize reducing the LQR cost by increasing $\epsilon$, the blade increases the outboard chord and the twist, as shown in \cref{fig:quadrotor-rotor-design-3D,fig:quadrotor-rotor-design}.
The chord reaches the upper bound of 35\% radius in the inboard section, and therefore, we do not see the variation in the inboard chord.
The larger chord and twist allow the rotors to generate higher thrust under the same control inputs, enabling quicker attitude recovery.
This lowers the LQR cost by 10\% while compromising on the hovering power by 3\%.

\begin{figure}[htbp]
    \centering
    \includegraphics[width=1.0\columnwidth]{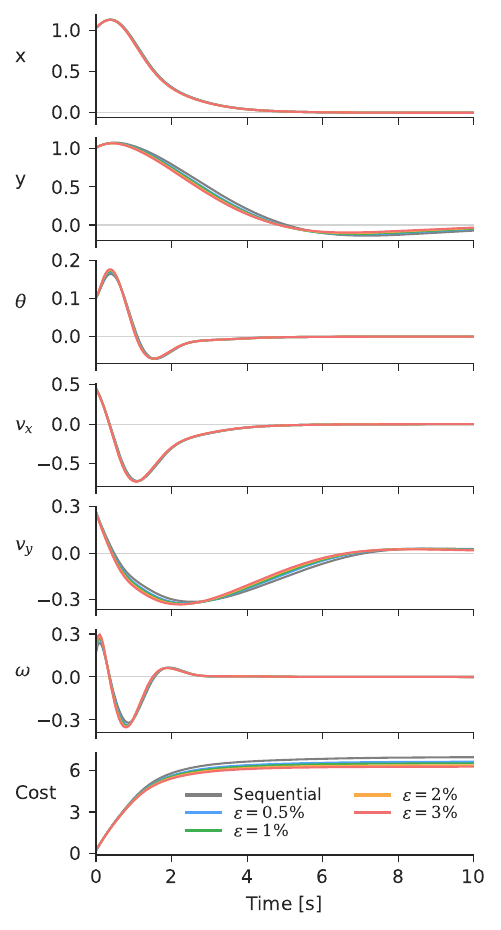}
    \caption{State time histories for the quadrotor attitude recovery, where each line corresponds to a point in \cref{fig:quadrotor-pareto-front} with the same color.}
    \label{fig:quadrotor-state-history}
\end{figure}

\Cref{fig:quadrotor-state-history} shows the time histories of the state variables, starting from the initial perturbed states to recover to the target zero states.
The state response becomes faster as we increase the $\epsilon$ value and tailor optimization toward reducing the LQR cost.
This trend is particularly visible in $y$ and $v_y$ histories in \cref{fig:quadrotor-state-history}.

Overall, the results of this quadrotor CCD problem demonstrates that CCD optimization is useful in investigating the trade-off between steady and dynamic performance, allowing us to make a more informed system-level design decision.
The adjoint-based sensitivity analysis we developed is essential in enabling such optimization studies within reasonable computational cost.

\section{Conclusion}
\label{sec:conclusion}

Control co-design optimization with gradient-based methods has been intractable for large-scale problems due to the prohibitive cost of computing derivatives with respect to large numbers of design variables.
We develop an efficient adjoint-based method that addresses this fundamental challenge by computing closed-loop system performance derivatives in a single pass, enabling practical optimization of systems with hundreds of design variables.

The main insight is that the coupled adjoint equations exhibit a special feed-forward structure that can be exploited to avoid solving a large monolithic system.
Rather than directly computing a coupled Jacobian, we decompose the problem into three smaller adjoint equations that can be solved sequentially using block back-substitution.
This approach eliminates the need for specialized coupled solvers with custom preconditioners, significantly reducing implementation complexity while maintaining computational efficiency.
Furthermore, one of these adjoint equations reduces to a Lyapunov equation, allowing us to leverage existing, well-established numerical infrastructure.

We demonstrate the practical effectiveness of this method on two representative problems. For the cart-pole system, we achieve a 76.79\% reduction in control cost through simultaneous design and control optimization.
For the quadrotor blade co-design problem, we reduce dynamic performance (LQR cost) by 10\% with only a 3\% penalty in steady-state hovering power—a trade-off that would be difficult to achieve without gradient-based optimization.

The framework is general and extends naturally to other closed-loop control formulations beyond LQR, such as model predictive control (MPC) and $\mathcal{H}_\infty$ control, opening opportunities for future research on gradient-based co-design of more sophisticated control strategies in applications including aircraft, wind turbines, and robotic systems.  


\section*{Acknowledgments}
The authors gratefully acknowledge Eirikur Jonsson for useful discussions.


\bibliographystyle{asmejour}   
\bibliography{references,mdolab}

\appendix

\section{Derivation of the adjoint equations using the Lagrangian} 
\label{ap:adjoint}

For a general function of interest $\hat{I}(\hat{\mb{d}},\hat{\mb{u}})$ and constraints $\hat{\mb{r}}(\hat{\mb{d}},\hat{\mb{u}})$, where $\hat{\mb{d}},\hat{\mb{u}}$ are design variable and state variable, respectively.
The Lagrangian can be written as
\begin{equation}
\hat{L} = \hat{I} + \hat{\boldsymbol{\psi}}^\intercal \hat{\mb{r}}.
\end{equation}
Then, we have that
\begin{equation}
\f{\d \hat{L}}{\d \hat{\mb{d}}} = \f{\d \hat{I}}{\d \hat{\mb{d}}} + \hat{\boldsymbol{\psi}}^\intercal \f{\d\hat{\mb{r}}}{\d \hat{\mb{d}}} = \f{\d \hat{I}}{\d \hat{\mb{d}}} = \f{\p \hat{I}}{\p \hat{\mb{d}}} + \f{\p \hat{I}}{\p \hat{\mb{u}}}\f{\d \mb{u}}{\d \mb{d}} + \hat{\boldsymbol{\psi}}^\intercal \f{\p\hat{\mb{r}}}{\p \hat{\mb{d}}} + \hat{\boldsymbol{\psi}}^\intercal \f{\p\hat{\mb{r}}}{\p \hat{\mb{u}}}\f{\d \hat{\mb{u}}}{\d \hat{\mb{d}}}.
\end{equation}
By regrouping the equation, we have
\begin{equation}
\f{\d \hat{I}}{\d \hat{\mb{d}}} = \left(\f{\p \hat{I}}{\p \hat{\mb{d}}} + \hat{\boldsymbol{\psi}}^\intercal \f{\p\hat{\mb{r}}}{\p \hat{\mb{d}}}\right) + \left(\f{\p \hat{I}}{\p \hat{\mb{u}}} + \hat{\boldsymbol{\psi}}^\intercal \f{\p\hat{\mb{r}}}{\p \hat{\mb{u}}}\right)\f{\d \hat{\mb{u}}}{\d \hat{\mb{d}}}.
\end{equation}
By enforcing the Lagrangian multiplier to satisfy the adjoint equation (it also makes Lagrange multiplier equal to the adjoint vector), we have
\begin{equation}
\begin{aligned}
\f{\d \hat{I}}{\d \hat{\mb{d}}} &= \f{\p \hat{I}}{\p \hat{\mb{d}}} + \hat{\boldsymbol{\psi}}^\intercal \f{\p\hat{\mb{r}}}{\p \hat{\mb{d}}},\\
0&=\f{\p \hat{I}}{\p \hat{\mb{u}}} + \hat{\boldsymbol{\psi}}^\intercal \f{\p\hat{\mb{r}}}{\p \hat{\mb{u}}},
\end{aligned}
\end{equation}
where the challenging to evaluet total derivative $\d \hat{\mb{u}} / \d \hat{\mb{d}}$ is cancelled out.

The adjoint equation can be obtained by setting $\p \hat{L} / \p \hat{\mb{u}} = 0$.
The following equation explicitly shows that
\begin{equation}
\f{\p \hat{L}}{\p \hat{\mb{u}}} = \f{\p \hat{I}}{\p \mb{u}} + \hat{\boldsymbol{\psi}}^\intercal \f{\p \hat{\mb{r}}}{\p \hat{\mb{u}}} = 0,
\end{equation}
which corresponds to the adjoint equation.

\section{Unsteady equation adjoint}
\label{sec:unsteady_adjoint}

The unsteady equation in forward mode can be written as
\begin{equation}
\begin{gathered}
{\mb{r}^{n_t}_{\mathrm{cl}}}^{(k)}=\delta\mb{x}^{(k)} - \delta\mb{x}^{(k-1)}\\ - \mb{r}_{\mathrm{nl}}(\mb{x}_{\mathrm{tgt}} + \delta\mb{x}^{(k-1)}, \mb{u}_{\mathrm{tgt}} + \mb{W}\delta\mb{x}^{(k-1)}, \mb{d})\Delta t = 0,\\ \,\, k = 1, \ldots, n,
\end{gathered}
\end{equation}
where the initial condition is set to
\begin{equation}
\delta\mb{x}^{(0)} = \delta\mb{x}_0 .
\end{equation}
Expressing this in matrix form we get
\begin{equation}
\begin{gathered}
\mb{r}_{\mathrm{cl}}^{n_t}
=
\begin{bmatrix}
\mb{I} & \mb{O} & \cdots & \mb{O} & \mb{O} \\
- \mb{I} & \mb{I} & \cdots & \mb{O} & \mb{O} \\
\vdots & \vdots & \ddots & \vdots & \vdots \\
\mb{O} & \mb{O} & \cdots & \mb{I} & \mb{O} \\
\mb{O} & \mb{O} & \cdots & - \mb{I} & \mb{I}
\end{bmatrix}
\begin{bmatrix}
\delta\mb{x}^{(1)}\\
\delta\mb{x}^{(2)}\\
\vdots\\
\delta\mb{x}^{(n-1)}\\
\delta\mb{x}^{(n)}
\end{bmatrix}
+
\\
\begin{bmatrix}
-\mb{I} \delta\mb{x}^{(0)} - \mb{r}_{\mathrm{nl}}(\mb{x}_{\mathrm{tgt}} + \delta\mb{x}^{(0)}, \mb{u}_{\mathrm{tgt}} + \mb{W}\delta\mb{x}^{(0)}, \mb{d})\Delta t \\
-\mb{r}_{\mathrm{nl}}(\mb{x}_{\mathrm{tgt}} + \delta\mb{x}^{(1)}, \mb{u}_{\mathrm{tgt}} + \mb{W}\delta\mb{x}^{(1)}, \mb{d})\Delta t \\
\vdots \\
-\mb{r}_{\mathrm{nl}}(\mb{x}_{\mathrm{tgt}} + \delta\mb{x}^{(n-2)}, \mb{u}_{\mathrm{tgt}} + \mb{W}\delta\mb{x}^{(n-2)}, \mb{d})\Delta t \\
-\mb{r}_{\mathrm{nl}}(\mb{x}_{\mathrm{tgt}} + \delta\mb{x}^{(n-1)}, \mb{u}_{\mathrm{tgt}} + \mb{W}\delta\mb{x}^{(n-1)}, \mb{d})\Delta t \\
\end{bmatrix}
=
\mb{0}.
\end{gathered}
\end{equation}
Thus, the adjoint equation for these unsteady time-accurate equations is
\begin{equation}
\f{\p \mb{r}_{\mathrm{cl}}^{n_t}}{\p \delta\mb{x}^{n_t}}^{\intercal}\boldsymbol{\psi}_{\mathrm{cl}} = -\f{\p I}{\p \delta\mb{x}^{n_t}}.
\end{equation}
In matrix form, we have
{\scriptsize
\begin{equation}
\begin{gathered}
\begin{bmatrix}
\mb{I} & - \mb{I} + \left(-{\mb{J}^{(1)}}^\intercal -\mb{W}^\intercal \mb{G^{(1)}}^\intercal\right)\Delta t & \cdots & \mb{O} & \mb{O} \\
\mb{O} & \mb{I} & \cdots & \mb{O} & \mb{O} \\
\vdots & \vdots & \ddots & \vdots & \vdots \\
\mb{O} & \mb{O} & \cdots & \mb{I} & - \mb{I} + \left(-{\mb{J}^{(n-1)}}^\intercal -\mb{W}^\intercal \mb{G^{(n-1)}}^\intercal\right)\Delta t \\
\mb{O} & \mb{O} & \cdots & \mb{O} & \mb{I}
\end{bmatrix}
\begin{bmatrix}
\boldsymbol{\psi}_{\mathrm{cl}}^{(1)}\\
\boldsymbol{\psi}_{\mathrm{cl}}^{(2)}\\
\vdots\\
\boldsymbol{\psi}_{\mathrm{cl}}^{(n-1)}\\
\boldsymbol{\psi}_{\mathrm{cl}}^{(n)}
\end{bmatrix}\\
=
\begin{bmatrix}
-\f{\p I}{\p \delta\mb{x}^{(1)}}\\
-\f{\p I}{\p \delta\mb{x}^{(2)}}\\
\vdots\\
-\f{\p I}{\p \delta\mb{x}^{(n-1)}}\\
-\f{\p I}{\p \delta\mb{x}^{(n)}}
\end{bmatrix},
\end{gathered}
\end{equation}}
where $\mb{J}^{(i)}, \mb{G}^{i}$ are defined by
{\scriptsize
\begin{equation}
\begin{aligned}
\mb{J}^{(i)} &\coloneqq \mb{J}(\mb{x}_{\mathrm{tgt}} + \delta\mb{x}^{(i)}, \mb{u}_{\mathrm{tgt}}+ \mb{W}\delta\mb{x}^{(i)}, \mb{d}),\\
\mb{G}^{(i)} &\coloneqq \mb{G}(\mb{x}_{\mathrm{tgt}} + \delta\mb{x}^{(i)}, \mb{u}_{\mathrm{tgt}}+ \mb{W}\delta\mb{x}^{(i)}, \mb{d}),
\end{aligned}
\end{equation}}
with $i = 1, \ldots, n - 1$.
Using back-substitution, we can solve this equation starting from the last row and proceed backwards until we reach the first row.

\section{$\left({\p \mb{r}_{\mathrm{cl}}^{n_t}}/{\p \mb{P}}\right)^\intercal {\psi}_{\mathrm{cl}}^{n_t}$}
\label{sec:pr_pP_1}
Using a forward Euler scheme, the term $\left({\p \mb{r}_{\mathrm{cl}}^{n_t}}/{\p \mb{P}}\right)^\intercal\boldsymbol{\psi}_{\mathrm{cl}}^{n_t}$ can be expanded as
\begin{equation}
\f{\p \mb{r}_{\mathrm{cl}}^{n_t}}{\p \mb{P}}^\intercal\boldsymbol{\psi}^{n_t}_{\mathrm{cl}} = \mb{G}_{\mathrm{tgt}}\mb{S}^{-\intercal}\sum_{k=1}^{n}{\mb{G}^{(k-1)}}^\intercal\boldsymbol{\psi}_{\mathrm{cl}}^{(k)}{\delta\mb{x}^{(k-1)}}^\intercal \Delta t,
\end{equation}
where the nonlinear dynamics is defined by \cref{eq:unsteady_control}.

\section{${\psi}_{\mathrm{cl}}^\intercal{\p\mb{r}_{\mathrm{cl}}^{n_t}}/{\p {\theta}}$, ${\p\mathrm{tr}\left({\Psi}_{\mathrm{ARE}}^\intercal\mb{R}_{\mathrm{ARE}}\right)}/{\p {\theta}}$}
\label{sec:pr_pP_2}
Using a forward Euler scheme, $\boldsymbol{\psi}_{\mathrm{cl}}^\intercal{\p\mb{r}_{\mathrm{cl}}^{n_t}}/{\p \boldsymbol{\theta}}$ can be computed using $\boldsymbol{\psi}_{\mathrm{cl}}^\intercal{\p\mb{r}_{\mathrm{cl}}^{n_t}}/{\p \mb{x}_{\mathrm{tgt}}}$ and $\boldsymbol{\psi}_{\mathrm{cl}}^\intercal{\p\mb{r}_{\mathrm{cl}}^{n_t}}/{\p \mb{u}_{\mathrm{tgt}}}$ shown as follows:
\begin{equation}
\begin{aligned}
\f{\p\mb{r}_{\mathrm{cl}}^{n_t}}{\p \mb{x}_{\mathrm{tgt}}}^\intercal\boldsymbol{\psi}^{n_t}_{\mathrm{cl}} &= \sum_{k=1}^n \Delta t {\mb{J}^{(i-1)}}^\intercal\boldsymbol{\psi}^{(i)}_{\mathrm{cl}}, \\
\f{\p\mb{r}_{\mathrm{cl}}^{n_t}}{\p \mb{u}_{\mathrm{tgt}}}^\intercal\boldsymbol{\psi}^{n_t}_{\mathrm{cl}} &= \sum_{k=1}^n \Delta t {\mb{G}^{(i-1)}}^\intercal\boldsymbol{\psi}^{(i)}_{\mathrm{cl}}. \\
\end{aligned}
\end{equation}

Similarly, the term ${\p\mathrm{tr}\left(\boldsymbol{\Psi}_{\mathrm{ARE}}^\intercal\mb{R}_{\mathrm{ARE}}\right)}/{\p \boldsymbol{\theta}}$ can be computed using ${\p\mathrm{tr}\left(\boldsymbol{\Psi}_{\mathrm{ARE}}^\intercal\mb{R}_{\mathrm{ARE}}\right)}/{\p \mb{x}_{\mathrm{tgt}}}$ and ${\p\mathrm{tr}\left(\boldsymbol{\Psi}_{\mathrm{ARE}}^\intercal\mb{R}_{\mathrm{ARE}}\right)}/{\p \mb{u}_{\mathrm{tgt}}}$
\begin{equation}
\begin{aligned}
\f{\p\mathrm{tr}\left(\boldsymbol{\Psi}_{\mathrm{ARE}}^\intercal\mb{R}_{\mathrm{ARE}}\right)}{\p \mb{x}_{\mathrm{tgt}}} &= \f{\p \mathrm{tr}\left(\mb{J}_{\mathrm{tgt}}^\intercal\left(\mb{P}\boldsymbol{\Psi}_{\mathrm{ARE}}^\intercal + \mb{P}^\intercal \boldsymbol{\Psi}_{\mathrm{ARE}}\right) - \mb{P}\mb{G}_{\mathrm{tgt}}\mb{S}^{-1}\mb{G}_{\mathrm{tgt}}^\intercal \mb{P}\right)}{\p \mb{x}_{\mathrm{tgt}}},\\
\f{\p\mathrm{tr}\left(\boldsymbol{\Psi}_{\mathrm{ARE}}^\intercal\mb{R}_{\mathrm{ARE}}\right)}{\p \mb{u}_{\mathrm{tgt}}} &= \f{\p \mathrm{tr}\left(\mb{J}_{\mathrm{tgt}}^\intercal\left(\mb{P}\boldsymbol{\Psi}_{\mathrm{ARE}}^\intercal + \mb{P}^\intercal \boldsymbol{\Psi}_{\mathrm{ARE}}\right) - \mb{P}\mb{G}_{\mathrm{tgt}}\mb{S}^{-1}\mb{G}_{\mathrm{tgt}}^\intercal \mb{P}\right)}{\p \mb{u}_{\mathrm{tgt}}}.\\
\end{aligned}
\end{equation}

\section{${\p\mb{r}_{\mathrm{cl}}^{n_t}}/{\p \mb{d}}^\intercal{\psi}_{\mathrm{cl}}^{n_t}$}
\label{sec:pr_pP_3}
And the term ${\boldsymbol{\psi}_{\mathrm{cl}}^{n_t}}^\intercal\left({\p \mb{r}_{\mathrm{cl}}^{n_t}}/{\p \mb{d}}\right)$ can be further expanded as
\begin{equation}
\f{\p\mb{r}_{\mathrm{cl}}^{n_t}}{\p \mb{d}}^\intercal\boldsymbol{\psi}_{\mathrm{cl}}^{n_t} = 
\sum_{k=1}^n\left( - \Delta t \f{\p \mb{r}_{\mathrm{nl}}}{\p\mb{d}}^\intercal(\mb{x}_{\mathrm{tgt}} + \delta\mb{x}^{(k - 1)}, \mb{u}_{\mathrm{tgt}} + \mb{W}\delta \mb{x}^{(k-1)}, \mb{d})\boldsymbol{\psi}_{\mathrm{cl}}^{(i)} \right).
\end{equation}

\section{Partial derivatives of \cref{eq:cost_func}}
\label{sec:FoI_partials}
We show the partial derivatives of the cost function defined by \cref{eq:cost_func}.
We assume that a forward Euler scheme is used for the unsteady equation.
The partial derivatives are
\begin{equation}
\begin{aligned}
\f{\p f}{\p \delta\mb{x}^{(i)}} &= 2\left(\mb{Q} + \mb{P}^\intercal\mb{G}_{\mathrm{tgt}}(\mb{x}_{\mathrm{tgt}}, \mb{d})\mb{S}^{-\intercal}\mb{G}_{\mathrm{tgt}}^\intercal(\mb{x}_{\mathrm{tgt}}, \mb{d}) \mb{P} \right) \delta\mb{x}^{(i)}\Delta t, \\
\f{\p f}{\p \mb{P}} &= \sum_{i=1}^{n_t}\mb{G}_{\mathrm{tgt}}(\mb{S}^{-\intercal} + \mb{S}^{-1})\mb{G}_{\mathrm{tgt}}^{\intercal}\mb{P}\delta\mb{x}^{(i)}{\delta\mb{x}^{(i)}}^\intercal\Delta t,\\
\f{\p f}{\p \mb{G}_{\mathrm{tgt}}} &= \sum_{i=1}^{n_t}\mb{P} \delta\mb{x}^{(i)}{\delta\mb{x}^{(i)}}^\intercal \mb{P}^\intercal \mb{G}_{\mathrm{tgt}}(\mb{S}^{-1} + \mb{S}^{-\intercal})\Delta t.
\end{aligned}
\end{equation}


\end{document}